# Metal-insulator transition for the almost Mathieu operator

By Svetlana Ya. Jitomirskaya*

**Abstract**

We prove that for Diophantine $\omega$ and almost every $\theta$, the almost Mathieu operator, $(H_{\omega,\lambda,\theta}\Psi)(n) = \Psi(n+1) + \Psi(n-1) + \lambda\cos 2\pi(\omega n + \theta)\Psi(n)$, exhibits localization for $\lambda > 2$ and purely absolutely continuous spectrum for $\lambda < 2$. This completes the proof of (a correct version of) the Aubry-André conjecture.

## 1. Introduction

The almost Mathieu operator $H_{\omega,\lambda,\theta}$, acting on $\ell^2(\mathbb{Z})$ and given by

$$(1.1) \qquad (H_{\omega,\lambda,\theta}\Psi)(n) = \Psi(n+1) + \Psi(n-1) + f(\omega n + \theta)\Psi(n),$$

with $f(\theta) = \lambda\cos(2\pi\theta)$, was first introduced by Peierls [39] and extensively studied in physics and mathematics literature since the 1970's. In this paper we discuss the decomposition of spectral measures of $H_{\omega,\lambda,\theta}$. For background and some recent results on other interesting topics, not mentioned here, see [26], [37], [30], [40], [43], [32], [33], [23].

Our main result is the following:

THEOREM 1. *For almost every $\omega \in \mathbb{R}, \theta \in \mathbb{R}$, the almost Mathieu operator $H_{\omega,\lambda,\theta}$ has*

$1^o$ *For $\lambda > 2$, only pure point spectrum with exponentially decaying eigenfunctions,*

$2^o$ *For $\lambda = 2$, purely singular-continuous spectrum,*

$3^o$ *For $\lambda < 2$, purely absolutely continuous spectrum.*

*Remark.* Part $2^o$ is proved in [25] and stated here for completeness only.

*Alfred P. Sloan Research Fellow. The author was supported in part by NSF Grant DMS-9704130.

Parts $1^o$ and $3^o$ of Theorem 1 were claimed by Aubry-André [1] (see also [2]) to hold for all $\theta$ and all irrational $\omega$. It was soon realized by Avron-Simon [4], based on the Gordon lemma [24], that $1^o$ does not hold for Liouville $\omega$, for which the spectrum of $H_{\omega,\lambda,\theta}$ is purely singular continuous. It was later understood [35] that arithmetic properties of $\theta$ are also important, since for a dense $G_\delta$ of $\theta$ for any given irrational $\omega$, operator $H_{\omega,\lambda,\theta}$ has purely singular-continuous spectrum. Thus an "a.e. $\omega,\theta$" part is necessary for $1^o$. We believe however that $3^o$ holds for all $\omega,\theta$, and $2^o$ for all $\theta$ and all irrational $\omega$.

What Aubry-André actually proved for $\lambda > 2$, was positivity of Lyapunov exponents which, by the general Ishii-Pastur theorem (see, e.g., [11]) implies absence of absolutely continuous spectrum (for a.e. $\theta$, and, by a recent result of Last-Simon [38], for all $\theta$). Their proof was made rigorous by Avron-Simon [3] and Pastur-Figotin [17], and a different proof followed from Herman [27]. Using the Dinaburg-Sinai technique [13], Bellissard-Lima-Testard [6] proved the existence of some pure point spectrum for $\lambda$ very large (and of some absolutely continuous spectrum for $\lambda$ very small) for a.e. $\omega$. The first proofs of complete localization (pure point spectrum with exponentially decaying eigenfunctions) for $\lambda$ large (a.e. $\omega,\theta$) are due to Sinai [42] and Fröhlich-Spencer-Wittwer [20]. Later, alternative KAM-type arguments for $\lambda$ large were developed by Eliasson [16] and Goldshtein [22]. For $\lambda$ very small, pure absolutely continuous spectrum was obtained by Chulaevsky-Delyon [8] using duality and the construction of Sinai [42]. Alternative proofs were given by Elliasson [15] and Goldshtein [22]. All those results use perturbative arguments and most hold for larger classes of quasiperiodic potentials (1.1).

The first nonperturbative results for $\lambda < 2$ were obtained by Last [36], who proved the existence of a large absolutely continuous component for any $\omega$ and any $\lambda < 2$ (and all $\theta$, as shown by Gesztesy-Simon [21] or as follows from [38]). In [28] a nonperturbative approach to localization was developed and pure point spectrum was proved for $\lambda \geq 15$. In [25] Gordon et al. established the dual result: pure absolutely continuous spectrum for $\lambda \leq 4/15$. In [29] the technique of [28] was combined with duality and the result of [36] to prove the existence of a large pure point component for any $\lambda > 2$, a.e. $\omega,\theta$. However, neither [36] nor [29] ruled out the existence of singular-continuous spectrum. Thus, the absence of singular continuous spectrum for $4/15 < \lambda < 15$, $\lambda \neq 2$, remained an issue, and it is solved in the present paper.

We note that $3^o$ immediately follows from $1^o$ by the strong version of duality [25], which states that for any irrational $\omega$, if $H_{\omega,\frac{4}{\lambda},\theta}$ has only pure point spectrum for a.e. $\theta$, then $H_{\omega,\lambda,\theta}$ has purely absolutely continuous spectrum for a.e. $\theta$. Therefore, we will concentrate on the proof of $1^o$.

We introduce the arithmetic conditions on $\omega,\theta$ that will enable us to prove point spectrum. We will say that $\omega$ is Diophantine if there exist $c(\omega) > 0$ and

$1 < r(\omega) < \infty$ such that

$$(1.2) \qquad |\sin 2\pi j\omega| > \frac{c(\omega)}{|j|^{r(\omega)}}$$

for all $j \neq 0$. It is well known that a.e. $\omega$ is Diophantine. We define the set of resonant phases:

$$(1.3) \qquad \Theta = \{\theta : \text{ the relation } |\sin 2\pi(\theta + (k/2)\omega)| \\ < \exp(-k^{\frac{1}{2r(\omega)}}) \text{ holds for infinitely many } k\text{'s}\}.$$

It follows from the Borel-Cantelli lemma that $\Theta$ has zero Lebesgue measure. It is, however, a dense $G_\delta$. Now we can formulate a more detailed version of part $1^o$ of Theorem 1.

THEOREM 2. *Suppose $\omega$ is Diophantine, $\theta \notin \Theta$, $\lambda > 2$. Then $H_{\omega,\lambda,\theta}$ has only pure point spectrum with exponentially decaying eigenfunctions.*

Moreover, it is possible to show that the eigenfunctions decay at exactly the Lyapunov rate. With the Lyapunov exponent $\gamma(E)$ defined by (2.2), we have the following:

THEOREM 3. *With $\omega, \lambda, \theta$ as in Theorem 2, for any eigenvalue $E$ of operator $H_{\omega,\lambda,\theta}$ the corresponding eigenfunction $\Psi_E$ satisfies*

$$(1.4) \qquad \lim_{|n|\to\infty} \frac{\ln(\Psi_E^2(n) + \Psi_E^2(n+1))}{2|n|} = -\gamma(E).$$

*Remarks.* 1. As discussed above, Theorem 2 is not true for Liouville $\omega$, nor does it hold for every $\theta$ for a given Diophantine $\omega$. Furthermore, it is possible to obtain a complete result on localization that is true for all irrational $\omega$ (new critical constants where spectrum changes from singular-continuous to pure point will appear, depending on the exponential rate of approximation of $\omega$ by rationals [31]).

2. The set $\Theta$, although quite small, is not optimal here. One can make it smaller by replacing $e^{-k^{\frac{1}{2r}}}$ in (1.3) by $e^{-k\varepsilon}$ for small $\varepsilon$; however that requires a more delicate argument. On the other hand, replacing it with $e^{-k\varepsilon}$ for sufficiently large $\varepsilon$, will lead to singular-continuous spectrum [35]. It was conjectured in [30] that there is a sharp threshold in $\varepsilon$ for the transition from pure point to singular-continuous spectrum. In any case, we did not try to optimize $\Theta$ in the present proof.

3. We prove Theorems 2 and 3 for any $\lambda$ such that the Lyapunov exponent $\gamma(E) > 0$ for all $E$. We believe that under such condition (and $\omega, \theta$ as in Theorem 2) this result should be true for operators (1.1) with rather general quasiperiodic potentials $f(\theta)$, just as positivity of the Lyapunov exponent implies zerodimensionality of spectral measures for all $\omega, \theta$ [34]. While positivity of the Lyapunov exponents can easily be extracted from the present proof, we only use it, rather than prove it. Nor do we use the $\lambda > 2$ condition otherwise.

In Section 2 we prove Theorems 2 and 3 up to the main technical Lemma 4. We prove Lemma 4 in Section 3 based on two other statements, Lemmas 5 and 7. The latter lemmas are proved, correspondingly, in Sections 4 and 5. In Sections 4 and 5 we also formulate and prove some lemmas that are more general than needed for the present proof. We believe they may be useful in certain other situations.

The author is grateful to the referee for the suggestions that led to significant improvement in this paper.

## 2. General setup. Proof of Theorems 2 and 3

*Definition.* A formal solution $\Psi_E(x)$ of the equation $H_{\omega,\lambda,\theta}\Psi_E = E\Psi_E$ will be called *a generalized eigenfunction* if $|\Psi_E(x)| \leq C(1 + |x|)$ for some $C = C(\Psi_E) < \infty$.

It is well-known that to prove Theorem 2 it suffices to prove that generalized eigenfunctions decay exponentially [7], [41].

We will use the notation $G_{[x_1,x_2]}(E)(x,y)$ for the Green's function $(H - E)^{-1}(x,y)$ of the operator $H_{\omega,\lambda,\theta}$ restricted to the interval $[x_1, x_2]$ with zero boundary conditions at $x_1 - 1$ and $x_2 + 1$. We now fix $E, \lambda \in \mathbb{R}$ and $\omega$ satisfying (1.2). To simplify notation, in some cases the $E, \lambda, \omega$-dependence of various quantities will be omitted. We denote the one-step transfer-matrix of $H\Psi = E\Psi$ by

$$B(\theta, E) = \begin{pmatrix} E - \lambda \cos 2\pi\theta & -1 \\ 1 & 0 \end{pmatrix}.$$

The $k$-step transfer-matrix is given by

$$M_k(\theta, E) = B(\theta + (k-1)\omega, E) \cdots B(\theta + \omega, E)B(\theta, E).$$

Let us denote

$$P_k(\theta, E) = \det\left[(H_{\omega,\lambda,\theta} - E)\Big|_{[0,k-1]}\right].$$

Then the $k$-step transfer-matrix can be written as

$$(2.1) \qquad M_k(\theta, E) = \begin{pmatrix} P_k(\theta, E) & -P_{k-1}(\theta + \omega, E) \\ P_{k-1}(\theta, E) & -P_{k-2}(\theta + \omega, E) \end{pmatrix}.$$

The Lyapunov exponent $\gamma(E)$ is defined as

$$(2.2) \qquad \gamma(E) = \lim_{n \to \infty} \frac{\int_0^1 \ln \|M_k(\theta, E)\| d\theta}{k} = \inf_k \frac{\int_0^1 \ln \|M_k(\theta, E)\| d\theta}{k}.$$

The limit exists and the last equality holds by the subadditive ergodic theorem (see, e.g., [11]).

Let $K = \{k \in \mathbb{N} : \text{ there exists } \theta \in [0,1] \text{ with } |P_k(\theta, E)| \geq \frac{1}{\sqrt{2}} e^{k\gamma(E)}\}$. By (2.1) for every $k \in \mathbb{N}$ we have that at least one of $k, k+1, k+2$ belongs to $K$.

*Definition.* Fix $E \in \mathbb{R}$, $m \in \mathbb{R}$. A point $y \in \mathbb{Z}$ will be called $(m,k)$-regular if there exists an interval $[x_1, x_2]$, $x_2 = x_1 + k - 1$, containing $y$, such that

$$|G_{[x_1,x_2]}(y, x_i)| < e^{-m|y-x_i|}, \text{ and } \text{dist}(y, x_i) \geq \frac{1}{5}k; \; i = 1, 2.$$

Otherwise, $y$ will be called $(m,k)$-singular.

This definition will be useful only for positive $m$; however, it can be used formally with negative $m$ as well. The number $\frac{1}{5}$ here can be replaced, equivalently, by any other number smaller than $\frac{1}{4}$. The main difference between this definition and that of regularity typically used in multiscale analysis is the flexibility allowed in the choice of the interval $[x_1, x_2]$. It is going to be essential for the proof.

It is well-known and can be checked easily that values of any formal solution $\Psi$ of the equation $H\Psi = E\Psi$ at a point $x \in I \subset \mathbb{Z}$ can be reconstructed from the boundary values via

$$(2.3) \qquad \Psi(x) = G_I(x, x_1)\Psi(x_1 - 1) + G_I(x, x_2)\Psi(x_2 + 1),$$

where $I = [x_1, x_2]$. This implies that if $\Psi_E$ is a generalized eigenfunction, then every point $y \in \mathbb{Z}$ with $\Psi_E(y) \neq 0$ is $(m,k)$-singular for $k$ sufficiently large: $k > k_1(E, m, \theta, y)$.

Theorem 2 will now follow from our next result.

LEMMA 4. *Suppose $\theta \notin \Theta$, $\omega$ satisfies (1.2). Then for every $y \in \mathbb{Z}$, $\varepsilon > 0, \alpha < 2$, there exists $k_2(\theta, \omega, y, \varepsilon, \alpha, E)$, such that for all $k \in K$, $k > k_2(\theta, \omega, y, \varepsilon, \alpha, E)$, if $x$ and $y$ are both $(\gamma(E) - \varepsilon, k)$-singular and $|x - y| > \frac{k+1}{2}$, then $|x - y| > k^\alpha$.*

*Remarks.* 1. If one somewhat relaxes the nonresonant condition on the excluded set $\Theta$ (with $\Theta$ still having measure 0, namely replacing the subexponential function by polynomial decay in (1.3)) one can obtain exponential in $k$ repulsion of singular clusters in Lemma 4.

2. This lemma is close in spirit to the central lemmas of multi-scale analysis ([18]; see also [19], [20], [14]). However our proof of it does not involve multiple scales.

To complete the proof of Theorem 2, we let $E(\theta)$ be a generalized eigenvalue of $H_{\omega,\lambda,\theta}$, and $\Psi_E(x)$ the corresponding generalized eigenfunction. Assume $\Psi(0) \neq 0$. Then, by Lemma 4, if $|x| > \max(k_1(E, \gamma(E) - \varepsilon, \theta, 0),$ $k_2(\theta, \omega, 0, \varepsilon, 1.5, E)) + 1$, the point $x$ is $(\gamma(E) - \varepsilon, k)$-regular for some $k \in \{|x| - 1, |x|, |x| + 1\} \cap K \neq \emptyset$, since $0$ is $(\gamma(E) - \varepsilon, k)$-singular. We, therefore, obtain that there exists an interval $[x_1, x_2]$ containing $x$, such that

$$\text{(2.4)} \qquad \frac{1}{5}(|x| - 1) \leq |x_i - x| \leq \frac{4}{5}(|x| + 1),$$

$$|G_{[x_1, x_2]}(x, x_i)| \leq e^{-(\gamma(E) - \varepsilon)|x - x_i|}, \ i = 1, 2.$$

Using (2.3), we obtain the estimate:

$$|\Psi_E(x)| \leq 2C(\Psi_E)(2|x| + 1)e^{-(\frac{\gamma(E) - \varepsilon}{5})(|x| - 1)}.$$

Since, for $\lambda > 2$, we have $\gamma(E) > 0$ for all $E$ (note: this is the only place in the paper where we use the condition $\lambda > 2$), this implies exponential decay. □

*Proof of Theorem* 3. The fact that the lower limit of the expression in (1.4) is always not lower than $-\gamma(E)$ is a general statement, proved by Craig-Simon [10] for any quasiperiodic (with a.e. $\theta$ form true even for any ergodic) potential. Therefore, we are concerned here only with the upper bound.

It is one of the main technical points of multiscale analysis that the exponential decay of a Green's function at a scale $k$ under certain conditions generates exponential decay with the same rate at a scale $k^\alpha$. The proof is usually done using block-resolvent expansion, with the combinatorial factor being killed by the power-law growth of scales. The proof of Theorem 3 will consist, roughly, of adapting this type of argument to our situation. We would like to note, however, that power-law repulsion of singular clusters, as in Lemma 4, although used in the following proof, is not essential here. Namely, one can prove Theorem 3 using only the fact that $x$ is $(\gamma(E) - \varepsilon, |x|)$-regular for all sufficiently large $|x|$.

Fix $\alpha, \alpha_1$ so that $1 < \alpha < \alpha_1 < 2$. Find $\bar{k}$ such that for $k > \bar{k}$ we have $(k+3)^\alpha + 4k^\alpha + 5k^{\alpha - 1} < k^{\alpha_1}$. Assume without loss of generality that $x$ is positive. Let $E$ be an eigenvalue of $H_{\omega, \lambda, \theta}$. Since $\Psi_E \in \ell^2$ (it even decays exponentially, as we already proved), we can normalize it to have an *a priori* bound

$$\text{(2.5)} \qquad |\Psi_E(x)| \leq 1, \ x \in \mathbb{Z}.$$

Take $x$ large enough: $|x| > k_0^\alpha$, where

$$k_0 = \max(k_1(E, \gamma(E) - \varepsilon, \theta, 0), k_2(\theta, \omega, 0, \varepsilon, \alpha_1, E), \bar{k}) + 3.$$

Let $k \in K$ be such that

$$\text{(2.6)} \qquad k^\alpha < |x| \leq (k+3)^\alpha.$$

Then $k > k_0 - 3$, and we have, by Lemma 4, that every $y \in [\frac{k+3}{2}, k^{\alpha_1}]$ is $(\gamma(E)-\varepsilon, k)$-regular. For each such $y$ pick an interval $[x_1, x_2]$ from the definition of $(\gamma(E) - \varepsilon, k)$-regularity and denote it by $I(y)$. We denote the boundary of the interval $I(y)$, the set $\{x_1, x_2\}$, by $\partial I(y)$. For $z \in \partial I(y)$ we let $z'$ be the neighbor of $z$ (i.e., $|z - z'| = 1$) not belonging to $I(y)$.

We now expand $\Psi(x_1 - 1)$ in (2.3) iterating (2.3) with $I = I(x_1 - 1)$, and $\Psi(x_2 + 1)$ using (2.3) with $I = I(x_2 + 1)$. Note that $I(x_i \pm 1)$, $i = 1, 2$, are well-defined, since, by construction, $\frac{k+3}{2} < x_1 - 1 < x_2 + 1 < k^{\alpha_1}$, and therefore points of the form $x_i \pm 1$, $i = 1, 2$, are $(\gamma(E) - \varepsilon, k)$-regular. We continue to expand each term of the form $\Psi(z)$ in the same fashion until we arrive at such a $z$ that either $z < k$ or the number of $G_I$ terms in the product becomes $5k^{\alpha-1}$, whichever comes first. We then obtain an expression of the form
(2.7)
$$\Psi(x) = \sum_{s; z_{i+1} \in \partial I(z_i')} G_{I(x)}(x, z_1) G_{I(z_1')}(z_1', z_2) \cdots G_{I(z_s')}(z_s', z_{s+1}) \Psi(z_{s+1}'),$$
where in each term of the summation we have $z_i > k$, $i = 1, \ldots, s$, and either $0 < z_{s+1}' < k$, $s \leq 5k^{\alpha-1}$, or $s + 1 = 5k^{\alpha-1}$. Note, that since we bound the number of expansions by $5k^{\alpha-1}$ and since, by construction, $k > \bar{k}$, we have $\frac{k+3}{2} < k < z_i' < x + 4k^\alpha + 5k^{\alpha-1} < k^{\alpha_1}$ for each $i \leq s$. Therefore each $z_i'$, $i \leq s$, is $(\gamma(E) - \varepsilon, k)$-regular, and $I(z_i')$ is well-defined. We now consider the two cases, $0 < z_{s+1}' < k$, and $s + 1 = 5k^{\alpha-1}$, separately.

If $0 < z_{s+1}' < k$, we have, by the definition of regularity and (2.5),
$$|G_{I(x)}(x, z_1) G_{I(z_1')}(z_1', z_2) \cdots G_{I(z_s')}(z_s', z_{s+1}) \Psi(z_{s+1}')|$$
$$\leq e^{-(\gamma(E)-\varepsilon)(|x-z_1| + \sum_{i=1}^{s} |z_i' - z_{i+1}|)}$$
$$\leq e^{-(\gamma(E)-\varepsilon)(|x - z_{s+1}| - (s+1))} \leq e^{-(\gamma(E)-\varepsilon)(x - k - 5k^{\alpha-1})}.$$

If $s + 1 = 5k^{\alpha-1}$, again by (2.5), the definition of regularity, and the fact that $|x - z_1| \geq \frac{k}{5}$, $|z_i - z_{i+1}| \geq \frac{k}{5}$, $i = 1, \ldots, s$, we can estimate
$$|G_{I(x)}(x, z_1) G_{I(z_1')}(z_1', z_2) \cdots G_{I(z_s')}(z_s', z_{s+1}) \Psi(z_{s+1}')| \leq e^{-(\gamma(E)-\varepsilon)\frac{k}{5} 5k^{\alpha-1}}.$$

Using (2.6) we obtain that in either case
(2.8)
$$|G_{I(x)}(x, z_1) G_{I(z_1')}(z_1', z_2) \cdots G_{I(z_s')}(z_s', z_{s+1}) \Psi(z_{s+1}')| \leq e^{-(\gamma(E)-\varepsilon-\delta)x}$$
for any $\delta > 0$ and $x$ sufficiently large. Finally, we observe that the total number of terms in (2.7) is bounded above by $2^{5k^{\alpha-1}}$. Combining it with (2.7), (2.8),(2.6), we obtain
$$|\Psi(x)| \leq 2^{5x^{\frac{\alpha-1}{\alpha}}} e^{-(\gamma(E)-\varepsilon-\delta)x}.$$
Since $\delta$ and $\epsilon$ can be chosen arbitrarily small for sufficiently large $x$, this implies the needed upper bound. □

## 3. Proof of Lemma 4

In this section we assume $E$ fixed and drop it from most of the notation.

We start with outlining the idea of the proof. Every $(m, k)$-singular point $x$ can be understood, very roughly, as a local center of localization, at the rate $m$ and on the scale $k$. If so, for any $x_1 < x < x_1 + k - 1$, such that both $x_1$ and $x_1 + k - 1$ are sufficiently far from $x$, the operator $H_{\omega,\lambda,\theta}$ restricted to the interval $[x_1, x_1 + k - 1]$ should have an eigenvalue close to $E$. This could be manifested by $|P_k(\theta + x_1\omega, E)|$ being abnormally small (exponentially smaller than $e^{\gamma(E)k}$). Then, for each singular point $x$ there exists a sufficiently long (longer than $\frac{k+1}{2}$) interval $I(x) \subset \mathbb{Z}$ of values of $x_1$ such that $|P_k(\theta + x_1\omega, E)|$ is small. This is made precise in Lemma 6. The main technical difficulty in the proof of this lemma lies in establishing a uniform (in $\theta$) upper bound on $|P_k(\theta, E)|$, that is a subject of Lemma 5, the proof of which is given in Section 4.

The rest of the proof consists mainly of establishing the following fact: if the points $\theta_1, \ldots, \theta_{k+1}$ are, in a certain sense, uniformly distributed on $[0, 1)$, and do not come abnormally close to one another, then $P_k(\theta)$ will not be small (in the sense above) at at least one of those points. Finally, for Diophantine $\omega$, any points of the form $\theta + x_1\omega$ with $x_1$ running through a union of two long enough intervals, $I(x)$ and $I(y)$, are going to be uniformly distributed (Lemma 12). If $x$ and $y$ are not too close to each other, the intervals $I(x)$ and $I(y)$ will not intersect and their union will consist of at least $k+1$ points. Also, the Diophantine condition on $\omega$ (1.2) together with the nonresonant condition on $\theta$ (1.3) imply that if $x$ is not very far from $y$, the points $\theta + x_1\omega$, with $x_1 \in I(x) \cup I(y)$, do not come too close to one another (Lemma 13). Summing it all up yields that two singular points, $x$ and $y$ (if they are not very close) must be sufficiently far from each other.

We will now proceed with a formal proof.

It is a standard linear algebra (Cramer's) rule that for any $x_1$, $x_2 = x_1 + k - 1$, $x_1 \leq y \leq x_2$,

$$
\begin{aligned}
(3.1) \qquad |G_{[x_1,x_2]}(x_1, y)| &= \left| \frac{P_{x_2-y}(\theta + (y+1)\omega)}{P_k(\theta + x_1\omega)} \right|, \\
|G_{[x_1,x_2]}(y, x_2)| &= \left| \frac{P_{y-x_1}(\theta + x_1\omega)}{P_k(\theta + x_1\omega)} \right|.
\end{aligned}
$$

We will need the following bound, uniform in $\theta$, on the numerator in (3.1):

LEMMA 5. *Suppose $H_{\omega,\theta}$ is given by (1.1) with $f(\theta) = \sum_{k=0}^{p} a_k \cos^k \theta$ a trigonometric polynomial, $\omega$ any irrational. Then for every $E \in \mathbb{R}$, $\varepsilon > 0$, there exists $k_f(\varepsilon, E, \omega)$ such that $|P_n(\theta)| < e^{(\gamma(E)+\varepsilon)n}$ for all $n > k_f(\varepsilon, E, \omega)$, all $\theta$.*

This lemma will be proved in Section 4.

Let $A_k^{z,\theta} = \{x \in \mathbb{Z} : |P_k(\theta + x\omega)| \leq e^{kz}\}$. We will show that every singular point "produces" a cluster of points belonging to an appropriate $A_k^z$; however, for $z < \gamma(E)$, such clusters in the set $A_k^z$ should be far apart.

LEMMA 6. *Suppose* $y \in \mathbb{Z}$ *is* $(\gamma(E) - \varepsilon, k)$-*singular,* $\varepsilon < \frac{\gamma(E)}{3}$, $k > 4k_f(\frac{\varepsilon}{6}, E, \omega) + 1$. *Then for any $x$ such that* $y - [\frac{3}{4}]k \leq x \leq y - [\frac{3}{4}]k + [\frac{k+1}{2}]$, *$x$ belongs to* $A_k^{\gamma(E) - \frac{\varepsilon}{8}, \theta}$.

*Proof.* This follows immediately from the definition of regularity, (3.1), and Lemma 5. □

Assume now that $y_1$ and $y_2$ are both $(\gamma(E) - \varepsilon, k)$-singular. Assume without loss of generality that $y_2 > y_1$. We set $d = y_2 - y_1$, $x_i = y_i - [3/4k]$, $i = 1, 2$. Note that $P_k(\theta)$ is an even function of $\theta + \frac{k-1}{2}\omega$ and can be written as a polynomial of degree $k$ in $\cos 2\pi(\theta + \frac{k-1}{2}\omega)$:

$$P_k(\theta) = \sum_{j=0}^{k} c_j \cos^j 2\pi(\theta + \frac{k-1}{2}\omega) \stackrel{\text{def}}{=} Q_k(\cos 2\pi(\theta + \frac{k-1}{2}\omega)).$$

We now set

$$\theta_j = \begin{cases} \theta + (x_1 + \frac{k-1}{2} + j)\omega, & j = 0, 1, \ldots, [\frac{k+1}{2}] - 1 \\ \theta + (x_2 + \frac{k-1}{2} + j - [\frac{k+1}{2}])\omega, & j = [\frac{k+1}{2}], \ldots, k. \end{cases}$$

By the assumption, $d > \frac{k+1}{2}$. This implies that all $\theta_j$, $j = 0, \ldots, k$, are different. We have, by Lemma 6, that $|Q_k(\cos 2\pi\theta_j)| < e^{k(\gamma(E) - \varepsilon/8)}$, $j = 0, \ldots, k$. We now write the polynomial $Q_k(z)$ in Lagrange interpolation form using $\cos 2\pi\theta_0, \ldots, \cos 2\pi\theta_k$:

$$(3.2) \qquad |Q_k(z)| = \left|\sum_{j=0}^{k} Q_k(\cos 2\pi\theta_j) \frac{\prod_{\ell \neq j}(z - \cos 2\pi\theta_\ell)}{\prod_{\ell \neq j}(\cos 2\pi\theta_j - \cos 2\pi\theta_\ell)}\right|.$$

The points $\theta_j$, representing two long pieces of the trajectory of an ergodic rotation, are uniformly distributed on the circle, and, as we will show in Lemma 13, the $\cos 2\pi\theta_j$ do not come too close to one another. We claim that this implies that the ratio of products in (3.2) does not deviate much from $e^0$. We make it precise in the following lemma that will be proved in Section 5.

LEMMA 7. *Suppose* $d < k^\alpha$ *for some* $\alpha < 2$, *and* $\theta, \omega$ *as in Theorem 2. Then for any* $\varepsilon > 0$ *there exists* $k_3(\varepsilon, \omega, \alpha, \theta, y_1)$ *such that for* $k > k_3(\varepsilon, \omega, \alpha, \theta, y_1)$, *for any $z \in [-1, 1]$,*

$$\frac{|\prod_{\ell \neq j}(z - \cos 2\pi\theta_\ell)|}{|\prod_{\ell \neq j}(\cos 2\pi\theta_j - \cos 2\pi\theta_\ell)|} \stackrel{\text{def}}{=} \frac{|I_1|}{|I_2|} \leq e^{k\varepsilon}.$$

We can now complete the proof of Lemma 4. Set $\hat{k} = \max\{4k_{\cos}(\frac{\varepsilon}{6}, E, \omega) + 1, k_3(\frac{\varepsilon}{16}, \omega, \alpha, \theta, y_1)\}$. Take $\alpha < 2$, and $k \in K$, $k > \hat{k}$, and $\theta_0$ such that $|P_k(\theta_0)| \geq \frac{1}{\sqrt{2}} \exp(k\gamma(E))$. Let $z_0 = \cos 2\pi(\theta_0 + \frac{k-1}{2}\omega)$. Suppose $d < k^\alpha$. We evaluate $Q_k(z_0)$ using (3.2) and Lemma 7 to estimate the ratio of products in (3.2). Now,

$$\frac{1}{\sqrt{2}} e^{k\gamma(E)} \leq |Q_k(z_0)| \leq e^{k(\gamma(E)-\frac{\varepsilon}{8})}(k+1)e^{\frac{k\varepsilon}{16}}$$

and the contradiction proves Lemma 4 with $k_2 = \max(\hat{k}, k_4(\varepsilon))$. □

## 4. Polynomials and measure. Proof of Lemma 5

We introduce the upper Lyapunov exponent, $\bar{\gamma}(E, \theta) = \overline{\lim}_{k \to \infty} \frac{\ln \|M_k(\theta, E)\|}{k}$. Craig-Simon [10] showed that for any ergodic potential $\bar{\gamma}(E, \theta) \leq \gamma(E)$ for a.e. $\theta$ (for all $\theta$ in the quasiperiodic case (1.1)). Thus, by (2.1), for every $\theta$ and sufficiently large $n$, we have $|P_n(\theta)| < e^{(\gamma(E)+\varepsilon)n}$. As before, $P_k(\theta)$ is an even function of $\theta + \frac{k-1}{2}\omega$ and can be written as a polynomial of degree $kp$ in $\cos 2\pi(\theta + \frac{k-1}{2}\omega)$. We will show that this implies that for large $n$, $|P_n(\theta)|$ can be bounded above uniformly in $\theta$. The key statement is:

THEOREM 8. *Let $Q(x) = \sum_{j=0}^{n} c_j x^j$ be an arbitrary $n^{\text{th}}$ degree polynomial. Suppose $|Q(z_0)| = a^n$, for some $a > 0$, $z_0 \in [-1, 1]$. Then for any $0 < b < a$ and sufficiently large $n$ ($n > n_0(a/b)$)*

$$|\{\theta \in (0, \pi) : |Q(\cos \theta)| < b^n\}| \leq c(a, b) < \pi.$$

Here and in the future we use $|\cdot|$ for Lebesgue measure. The nontrivial statement here is that the measure can be bounded by a constant not dependent on $n$ or the polynomial $Q$.

*Remark.* The theorem is certainly true for $|\{z \in [-1, 1] : |Q(z)| < b^n\}|$ as well. We formulate it with the cosine because it is convenient for our application and also simplifies the proof.

To complete the proof of Lemma 5 we set $A_n = \{\theta \in (0, \pi) : \text{for all } k > n, |P_k(\cos \theta)| \leq e^{(\gamma(E)+\varepsilon/2)k}\}$. Then $|A_n^c|$ goes to 0 as $n \to \infty$; however if there exists $\theta$ such that $|P_k(\cos \theta)| > e^{(\gamma(E)+\varepsilon)k}$, some $k > n$, then

$$\begin{aligned}|A_n^c| &\geq \left|\left\{\theta \in (0, \pi) : |P_k(\cos \theta)| > e^{(\gamma(E)+\varepsilon/2)k}\right\}\right| \\ &\geq \pi - c\left(e^{\frac{\gamma(E)+\varepsilon}{p}}, e^{\frac{\gamma(E)+\varepsilon/2}{p}}\right) > 0.\end{aligned}$$

This contradiction proves Lemma 5. □

*Proof of Theorem 8.* We start with the following elementary lemma, proven below:

LEMMA 9. *For any measurable set $B \subset \mathbb{R}$, $|B| = m$, any $n \in \mathbb{N}$, and any $\delta < m/n$ there exist $x_1, \ldots, x_{n+1} \in B$ such that $|x_i - x_j|/\delta \in \mathbb{N}$ for all $i \neq j$.*

In other words, one can find $n+1$ elements of an arithmetic progression with step $\delta$ belonging to $B$. We will now prove Theorem 8. Let $B = \{\theta \in (0, \pi) : |Q(\cos \theta)| < b^n\}$. Let $m = |B|$. We fix $\delta < \frac{m}{n}$, and pick, according to Lemma 9, $x_1, \ldots, x_{n+1} \in B$, elements of an arithmetic progression with step $\delta$. We write $Q(x)$ in a Lagrange interpolation form using $\cos x_1, \ldots, \cos x_{n+1}$. Then

$$(4.1) \quad a^n = |Q(z_0)| = \left| \sum_{j=1}^{n+1} Q(\cos x_j) \frac{\prod_{k \neq j}(z_0 - \cos x_k)}{\prod_{k \neq j}(\cos x_j - \cos x_k)} \right|$$

$$< (n+1)b^n \max_j \left| \frac{\prod_{k \neq j}(z_0 - \cos x_k)}{\prod_{k \neq j}(\cos x_j - \cos x_k)} \right|.$$

We denote the last ratio by $I$ and will estimate separately the numerator and denominator in $I$. First,

$$(4.2) \quad \exp\left(\sum_{j \neq k} \ln |z_0 - \cos x_j|\right) \leq \exp\left(\sum_{j=0}^{n} \ln |1 - \cos(\pi - j\delta)|\right)$$

$$\leq \exp\left(1/\delta \int_0^{n\delta} \ln(1 + \cos x)dx + \ln 2\right).$$

Similarly,

(4.3)

$$\exp\left(\sum_{j \neq k} \ln |\cos x_k - \cos x_j|\right) \geq \exp\left(\sum_{\substack{j = -[\frac{n+1}{2}] \\ j \neq 0}}^{[\frac{n+1}{2}]} \ln |\cos(\pi/2 + j\delta)|\right)$$

$$\geq \exp\left(2/\delta \int_0^{\frac{n+1}{2}\delta} \ln \sin x\, dx\right).$$

We set $m_1 = n\delta$, and, combining (4.2), and (4.3) we obtain

(4.4)
$$I \leq 2 \exp\left(\frac{n}{m_1}\left(-\int_{m_1}^{\pi} \ln(1 + \cos x)dx + 2\int_{m_1/2}^{\pi/2} \ln \sin x\, dx\right) - \ln \sin \frac{m_1}{2}\right)$$

where we used the elementary identity (5.8). We set

$$g(m) = \exp\left(\frac{1}{m}\left(-\int_m^{\pi} \ln(1 + \cos x)dx + 2\int_{m/2}^{\pi/2} \ln \sin x\, dx\right)\right).$$

It is easily verified that $g$ is a monotone decreasing function on $(0, \pi]$, and $g(\pi) = 1$. For $m_1 \geq \frac{\pi}{3}$ we have $\sin \frac{m_1}{2} \geq 1/2$ and we can rewrite (4.1) as

$$\frac{a}{b} \left( \frac{1}{4(n+1)} \right)^{\frac{1}{n}} < g(m_1)$$

which implies $m_1 \leq g^{-1}\left(\frac{a}{b}(\frac{1}{4(n+1)})^{\frac{1}{n}}\right)$ and, finally, for any $\varepsilon$ and sufficiently large $n$, we obtain $m \leq (1+\varepsilon)g^{-1}(\frac{a}{b})$. □

*Proof of Lemma* 9. We will need an elementary statement:

LEMMA 10. *Let* $A = \cup_{i=1}^{\infty} A_i$, *where the* $A_i$ *are measurable sets*, $i \in \mathbb{N}$; *suppose every* $x \in A$ *belongs to no more than* $k$ *different* $A_j$'*s. Then* $|A| \geq \frac{1}{k} \sum_{i=1}^{\infty} |A_i|$.

*Proof.* We use induction in $k$. Set $B_j = A_j \cap (\cup_{i=1}^{j-1} A_i)$, $j = 2, \ldots$, and write

$$\sum_{i=1}^{\infty} |A_i| = |A_1| + |A_2 \setminus A_1| + \cdots + |A_n \setminus \cup_{i=1}^{n-1} A_i| + \cdots + \sum_{i=2}^{\infty} |B_i|$$

$$= |A| + \sum_{i=2}^{\infty} |B_i| \leq |A| + (k-1)|A| = k|A|$$

since if $x \in \cap_{j=1}^{k} B_{i_j}$, $2 \leq i_1 < i_2 < \cdots < i_k$, then $x$ necessarily belongs to $\cap_{j=0}^{k} A_{i_j}$ for some $i_0 < i_1$, which is a contradiction. □

To complete the proof of Lemma 9 we let $A = [0, \delta]$ and set $A_i = (B - i\delta) \cap A$, $i \in \mathbb{Z}$. Assume the statement of the lemma is false, which means every $x \in A$ belongs to no more than $n$ different $A_j$'s. Then, by Lemma 10,

$$\frac{m}{n} > \delta = |A| \geq |\cup_{i=-\infty}^{\infty} A_i| \geq \frac{1}{n} \sum_{i=-\infty}^{\infty} |A_i| = \frac{1}{n} \sum_{i=-\infty}^{\infty} |B \cap (A + i\delta)| = \frac{m}{n}.$$

This contradiction proves the lemma. □

## 5. Uniform distribution and small denominators. Proof of Lemma 7

Let $\phi(n)$ be a monotone increasing function, $\phi(n) = o(n)$. We will say that $x_1, \ldots, x_n$ are $\phi$-uniformly distributed on $[a, b]$ if for any $f \in C([a, b])$ we have $|\sum_{i=1}^{n} f(x_i) - \frac{n}{b-a} \int_a^b f(x) dx| \leq \phi(n) \mathrm{Var} f$.

Clearly, if $x_1, \ldots, x_n$ are $\phi$-uniformly distributed and $y_1, \ldots, y_n$ are $\psi$-uniformly distributed, then $x_1, \ldots, x_n, y_1, \ldots, y_n$ are $\phi(n) + \psi(n)$-uniformly distributed. We will say that $f \in C([a, b])$ is normal if it takes each value no more than $p$ times (some $p = p(f)$), and the family $\{\ln(z - f(x))\}_{z \in [\min f, \max f]}$ is absolutely equicontinuous, by which we mean that for every $\varepsilon > 0$ there exists $\delta > 0$ such that

$$\left| \int_{|z-f(x)|<\delta} \ln|z - f(x)|dx \right| < \varepsilon/2 \text{ for all } z \in [\min f, \max f]. \tag{5.1}$$

It can be easily verified that every analytic function (in particular, $\cos 2\pi x$) is normal.

LEMMA 11. *Suppose $f(x)$ is normal, $\psi(n) = o(n), \phi(n) = o(n)$. Then for every $\varepsilon > 0$ and $n$ sufficiently large ($n > n_f(\varepsilon, \phi, \psi)$) for any $x_1, \ldots, x_n$ that are $\phi$-uniformly distributed on $[a, b]$ and for any $j_1, \ldots, j_{\psi(n)} \in [1, \ldots, n]$,*

$$\sum_{\substack{i=1 \\ i \neq j_1, \ldots, j_{\psi(n)}}}^{n} \ln|z - f(x_i)| \leq \frac{n}{b-a} \left( \int_a^b \ln|z - f(x)|dx + \varepsilon \right), \tag{5.2}$$

$$\sum_{\substack{i=1 \\ i \neq j_1, \ldots, j_{\psi(n)}}}^{n} \ln|z - f(x_i)| \geq \frac{n}{b-a} \left( \int_a^b \ln|z - f(x)|dx - \varepsilon \right)$$
$$+ \phi(n)p(f) \ln \min_{\substack{i=1,\ldots,n \\ i \neq j_1, \ldots, j_{\psi(n)}}} |z - f(x_i)|.$$

*Proof.* For $0 < \delta < 1$ put $g_\delta(z) = \begin{cases} \ln|z|, & |z| > \delta \\ \ln \delta, & |z| \leq \delta \end{cases}$. Given $\varepsilon > 0$, fix $\delta(\varepsilon)$ such that (5.1) is satisfied. Then, with $\delta = \delta(\varepsilon)$, we obtain:

$$\sum_{\substack{i=1 \\ i \neq j_1, \ldots, j_{\psi(n)}}}^{n} \ln|z - f(x_i)| \leq \sum_{\substack{i=1 \\ i \neq j_1, \ldots, j_{\psi(n)}}}^{n} g_\delta(|z - f(x_i)|)$$

$$\leq \frac{n}{b-a} \int_a^b g_\delta(|z - f(x)|)dx + \phi(n)\mathrm{Var}\, g_\delta(|z - f(x)|) - \psi(n)\ln\delta$$

$$\leq \frac{n}{b-a} \left( \int_a^b \ln(|z - f(x)|)dx + \varepsilon/2 \right)$$
$$+ \phi(n)p(f)(-\ln\delta + \ln(\max_x |z - f(x)|)) - \psi(n)\ln\delta$$

$$\leq \frac{n}{b-a} \left( \int_a^b \ln(|z - f(x)|)dx + \varepsilon \right), \text{ for } n > n_1(\varepsilon, \phi, \psi, f).$$

Let $\min_{\substack{i=1,\ldots,n \\ i \neq j_1, \ldots, j_{\psi(n)}}} |z - f(x_i)| = c < 1$. For an estimate from below we write:

$$\sum_{\substack{i=1 \\ i \neq j_1, \ldots, j_{\psi(n)}}}^{n} \ln |z - f(x_i)| = \sum_{\substack{i=1 \\ i \neq j_1, \ldots, j_{\psi(n)}}}^{n} g_c(|z - f(x_i)|)$$

$$\geq \frac{n}{b-a} \int_a^b g_c(|z - f(x)|)dx - \psi(n) \max_x \ln|z - f(x)|$$
$$- \phi(n) \mathrm{Var}_x g_c(|z - f(x)|)$$

$$\geq \frac{n}{b-a}\left(\int_a^b \ln|z - f(x)|dx - \varepsilon\right) + \phi(n)p(f)\ln c, \text{ for } n > n_2(\varepsilon, \phi, \psi). \quad \square$$

Let $\{x\}$ denote the fractional part of $x$. We will use the following:

LEMMA 12. *For any $\omega$ satisfying (1.2) and any $x \in \mathbb{R}$, points $\{x\}$, $\{x + \omega\}, \ldots, \{x + (n-1)\omega\}$ are $\phi$-uniformly distributed on $[0, 1]$, with $\phi(n) = c_1(\omega) n^{1 - r(\omega)^{-1}} \ln n$. Here $r(\omega)$ is given by (1.2).*

This is a slightly different form of a known statement (see, e.g., [33]). For the reader's convenience we provide a proof here:

*Proof.* Let $p_n/q_n$ be the sequence of continued fraction approximants of $\omega$. Let $n(k)$ be such that $q_{n(k)} \leq k < q_{n(k)+1}$. We will use $r$ for $r(\omega)$. Writing $k = b_n q_n + b_{n-1} q_{n-1} + \cdots + b_1 q_1 + b_0$ and using the Denjoy-Koksma inequality (see, e.g., Lemma 4.1, Ch. 3 of [9]), we get

$$(5.3) \quad \left|\sum_{j=0}^{k-1} f(\{\theta + j\omega\}) - k \int_0^1 f(\theta)d\theta\right| \leq (b_0 + \cdots + b_n)\mathrm{Var}(f)$$
$$\leq \left(\sum_{i=0}^{n}\left[\frac{q_{i+1}}{q_i}\right]\right)\mathrm{Var}(f).$$

Since (1.2) implies $q_{i+1} < \frac{q_i^r}{c}$, where $c = \frac{c(\omega)}{2\pi}$, we have $\frac{q_{i+1}}{q_i} < \frac{q_{i+1}}{(cq_{i+1})^{1/r}} = \frac{q_{i+1}^{1-1/r}}{c^{1/r}}$. The right-hand side of (5.3) can now be estimated as

$$\leq \left(c^{-1/r} \sum_{i=1}^{n(k)} q_i^{1-1/r} + \frac{k}{q_{n(k)}}\right)\mathrm{Var}(f) < \left(c^{-1/r} n(k) q_{n(k)}^{1-1/r} + \frac{k}{q_{n(k)}}\right)\mathrm{Var}(f).$$

Since $k < q_{n(k)+1} \leq \frac{q_{n(k)}^r}{c}$, we have $q_{n(k)} \geq (ck)^{1/r}$, and $\frac{k}{q_{n(k)}} \leq c^{-1/r} k^{1-1/r}$. Also, for any $\omega$, $n(k) \leq \frac{2 \ln q_{n(k)}}{\ln 2}$ (see, e.g. Ch. 4 of [9]). Thus we can continue our estimate as

$$\leq (c_1 k^{1-1/r} \ln k)\mathrm{Var}(f). \quad \square$$

Besides Lemma 12, the Diophantine properties of $\omega$ and the arithmetic properties of $\theta$ will only play a role in the following simple lemma that allows us to estimate the small denominators in (3.2).

LEMMA 13. *Suppose $\theta \notin \Theta$, $\omega$ satisfies (1.2); then for $k$ sufficiently large $(k > k_5(y_1, \theta, \omega))$,*

$$(5.4) \quad |\cos 2\pi\theta_j - \cos 2\pi\theta_\ell| \geq \exp\left(-(d+k/4+1)^{\frac{1}{2r(\omega)}}\right) \frac{c(\omega)}{(d+k+1)^{r(\omega)}}.$$

*Proof.* Define $i_j \in \mathbb{Z}$ by $i_j = \frac{\theta_j - \theta}{\omega} - (x_1 + \frac{k-1}{2})$, $j = 0, \ldots, k$. Then

$$(5.5) \quad i_j - x_1 \in [[-3k/4], \ldots, [-k/4]+1] \cap [d-[3k/4], \ldots, d-[k/4]+1].$$

We write
(5.6)
$$|\cos 2\pi\theta_j - \cos 2\pi\theta_\ell| = 2|\sin 2\pi(\theta + x_1\omega + \frac{k-1+i_j+i_\ell}{2}\omega) \sin 2\pi(i_j - i_\ell)\omega|.$$

Since $\theta \notin \Theta$ we have that for sufficiently large $j : j > k_6(\theta, \omega)$,

$$(5.7) \quad |\sin 2\pi(\theta + j/2)\omega| > \exp(-j^{\frac{1}{2r(\omega)}}).$$

Therefore we can combine (1.2), (5.7), (5.5), and (5.6) to obtain the statement of the lemma with $k_5(y_1, \theta, \omega) = k_6(\theta + y_1\omega, \omega)$. □

*Proof of Lemma 7.* Set $\phi(n) = 2c_1(\omega)(\frac{n}{2})^{1-r(\omega)^{-1}} \ln \frac{n}{2}$. According to Lemma 12 and a remark in the beginning of this section, points $\theta_0, \ldots, \theta_n$ are $\phi$-uniformly distributed on $[0, 1]$. Thus we can apply Lemma 11 with $f(x) = \cos x$, $\psi(n) = 1$, $\phi(n)$. The special property of the cosine that leads to spectral uniformity in energy (absence of mobility edge, in particular) is that

$$(5.8) \quad \int_0^1 \ln|z - \cos 2\pi x| dx = -\ln 2 \text{ for any } z \in [-1, 1].$$

The important part is that the integral in (5.8) is constant over $z$, not the actual value. In particular, it implies that $|I_1|$ and $|I_2|$ are of the same order. Precisely, we obtain, by Lemma 11, (5.8) and (5.4), that for $k$ large enough, $k > \max\left(n_{\cos}(\frac{\varepsilon}{3}, \phi(n), 1), k_5(y_1, \theta, \omega)\right)$, we have

$$(5.9) \quad |I_1| \leq \exp\left[(k+1)\left(-\ln 2 + \frac{\varepsilon}{3}\right)\right]$$

and

$$(5.10) \quad |I_2| \geq \exp\left[(k+1)\left(-\ln 2 - \frac{\varepsilon}{3}\right)\right.$$
$$\left. - 2\phi(k+1)\left((d+k/4+1)^{(2r(\omega))^{-1}} - \ln\frac{c(\omega)}{(d+k+1)^{r(\omega)}}\right)\right].$$

Suppose $d \leq k^\alpha$, some $\alpha < 2$. Then we can write $|I| = \frac{|I_1|}{|I_2|} \leq \exp((k+1)\frac{2\varepsilon}{3} + o(k))$ where $o(k)$ depends only on $\omega, \alpha$. This proves Lemma 7. □

UNIVERSITY OF CALIFORNIA, IRVINE, CALIFORNIA
*E-mail address*: szhitomi@uci.edu

## References

[1] S. AUBRY and G. ANDRÉ, Analyticity breaking and Anderson localization in incommensurate lattices, Ann. Israel Phys. Soc. **3** (1980), 133–164.
[2] S. AUBRY, The new concept of transitions by breaking of analyticity in a crystallographic model, Solid State Sci. **8** (1978), 264–277.
[3] J. AVRON AND B. SIMON, Almost periodic Schrödinger operators. II. The integrated density of states, Duke Math. J. **50** (1983), 369–391.
[4] ———, Singular continuous spectrum for a class of almost periodic Jacobi matrices, Bull. A.M.S. **6** (1982), 81–85.
[5] J. AVRON, P. VAN MOUCHE, and B. SIMON, On the measure of the spectrum for the almost Mathieu operator, Comm. Math. Phys. **132** (1990), 103–118; Erratum: ibid. **139** (1991), 215.
[6] J. BELLISARD, R. LIMA, and D. TESTARD, A metal-insulator transition for the almost Mathieu model, Comm. Math. Phys. **88** (1983), 207–234.
[7] Y. BEREZANSKII, *Expansions in Eigenfunctions of Selfadjoint Operators*, Transl. Math. Monogr., Vol. 17, A.M.S., Providence, RI, 1968.
[8] V. CHULAEVSKY and F. DELYON, Purely absolutely continuous spectrum for almost Mathieu operators, J. Statist. Phys. **55** (1989), 1279–1284.
[9] I. CORNFELD, S. FOMIN, and YA. SINAI, *Ergodic Theory*, Springer-Verlag, New York, 1982.
[10] W. CRAIG and B. SIMON, Subharmonicity of the Lyaponov index, Duke Math. J. **50** (1983), 551–560.
[11] H. L. CYCON, R. G. FROESE, W. KIRSCH, and B. SIMON, *Schrödinger Operators with Application to Quantum Mechanics and Global Geometry*, Springer-Verlag, New York, 1987.
[12] F. DELYON, Absence of localisation in the almost Mathieu equation, J. Phys. A **20** (1987), L21–L23.
[13] E. DINABURG and YA. SINAI, The one-dimensional Schrödinger equation with a quasiperiodic potential, Funct. Anal. Appl. **9** (1975), 279–289.
[14] H. VON DREIFUS and A. KLEIN, A new proof of localization in the Anderson tight binding model, Comm. Math. Phys. **124** (1989), 285–299.
[15] L. H. ELIASSON, Floquet solutions for the 1-dimensional quasiperiodic Schrödinger equation, Comm. Math. Phys. **146** (1992), 447–482.
[16] ———, Discrete one-dimensional quasi-periodic Schrödinger operators with pure point spectrum, Acta Math. **179** (1997), 153–196.
[17] A. FIGOTIN and L. PASTUR, The positivity of Lyapunov exponent and absence of the absolutely continuous spectrum for the almost-Mathieu equation, J. Math. Phys. **25** (1984), 774–777.
[18] J. FRÖHLICH and T. SPENCER, Absence of diffusion in the Anderson tight binding model for large disorder or low energy, Comm. Math. Phys. **88** (1983), 151–184.
[19] J. FRÖHLICH, F. MARTINELLI, E. SCOPPOLA, and T. SPENCER, Constructive proof of localization in the Anderson tight binding model, Comm. Math. Phys. **101** (1985), 21–46.
[20] J. FRÖHLICH, T. SPENCER, and P. WITTWER, Localization for a class of one-dimensional quasi-periodic Schrödinger operators, Comm. Math. Phys. **132** (1990), 5–25.
[21] F. GESZTESY and B. SIMON, The xi function, Acta Math. **176** (1996), 49–71.

[22] M. Goldshtein, Laplace transform method in perturbation theory of the spectrum of Schrödinger operator II (one-dimensional quasi-periodic potentials), preprint, 1992.
[23] F. Germinet, Dynamical localization II with an application to the almost Mathieu operator, J. Statist. Phys. **95**(1999), 273–286.
[24] A. Gordon, The point spectrum of the one-dimensional Schrödinger operator, Usp. Math. Nauk. **31** (1976), 257–258.
[25] A. Gordon, S. Jitomirskaya, Y. Last, and B. Simon, Duality and singular continuous spectrum in the almost Mathieu equation. Acta Math. **178** (1997), 169–183.
[26] B. Helffer and J. Sjöstrand, Semiclassical analysis for Harper's equation. III. Cantor structure of the spectrum, Mém. Soc. Math. France (N.S.) **39** (1989), 1–124.
[27] M. Herman, Une méthode pour minorer les exposants de Lyapounov et quelques exemples montrant le caractère local d'un théorème d'Arnol′d et de Moser sur le tore en dimension 2, Comment. Math. Helv. **58** (1983), 453–502.
[28] S. Jitomirskaya, Anderson localization for the almost Mathieu equation: a nonperturbative proof, Comm. Math. Phys. **165** (1994), 49–57.
[29] ______, Anderson localization for the almost Mathieu equation. II. Point spectrum for $\lambda > 2$, Comm. Math. Phys. **168** (1995), 563–570.
[30] ______, Almost everything about the almost Mathieu operator II, 373–382, *Proc.* XI *Internat. Congr. of Math. Phys.* (Paris, 1994), Internat. Press, Cambridge, MA, 1995.
[31] ______, Almost Mathieu equations with weakly Liouville frequencies: second critical constants, in preparation.
[32] S. Jitomirskaya and Y. Last, Dimensional Hausdorff properties of singular continuous spectra. Phys. Rev. Lett. **76** (1996), 1765–1769.
[33] ______, Anderson localization for the almost Mathieu equation, III. Semi-uniform localization, continuity of gaps, and measure of the spectrum, Comm. Math. Phys. **195** (1998), 1–14.
[34] ______, Power-law subordinacy and singular spectra, II. Line operators, preprint.
[35] S. Jitomirskaya and B. Simon, Operators with singular continuous spectrum, III. Almost periodic Schrödinger operators, Comm. Math. Phys. **165** (1994), 201–205.
[36] Y. Last, A relation between absolutely continuous spectrum of ergodic Jacobi matrices and the spectra of periodic approximants, Comm. Math. Phys. **151** (1993), 183–192.
[37] ______, Almost everything about the almost Mathieu operator I, 366–372, *Proc.* XI *Internat. Congr. Math. Phys.* (Paris, 1994), Internat. Press, Cambridge, MA, 1995.
[38] Y. Last and B. Simon, Eigenfunctions, transfer matrices, and absolutely continuous spectrum of one-dimensional Schrödinger operators, Invent. Math. **135** (1999), 329–367.
[39] R. Peierls, Zur Theorie des Diamagnetismus von Leitungselektronen, Z. Phys. **80** (1933), 763–791.
[40] M. Shubin, Discrete magnetic Laplacian, Comm. Math. Phys. **164** (1994), 259–275.
[41] B. Simon, Schrödinger semigroups, Bull. A.M.S. **7** (1982), 447–526.
[42] Ya. Sinai, Anderson localization for one-dimensional difference Schrödinger operator with quasiperiodic potential, J. Stat. Phys. **46** (1987), 861–909.
[43] P. B. Wiegmann, On the singular spectrum of the almost Mathieu operator – Arithmetics and Cantor spectra of integrable models, Prog. Theor. Phys. Suppl. **134** (1999), 171–181.